\documentclass[preprint,12pt]{elsarticle}
\graphicspath{ {./figures/} }
\usepackage{subfig}
\usepackage{xcolor}
\usepackage{amssymb}
\usepackage{amsmath}
\usepackage{hyperref}
\newtheorem{The}{Theorem}[section]
\newtheorem{Def}{Definition}[section]

\newtheorem{Ex}{Example}[section]
\journal{}
\begin{document}
	\begin{frontmatter}
		
		\title{Fractional Order Periodic Maps: Stability Analysis and Application to the Periodic-2 Limit Cycles in the Nonlinear Systems}
		
		\author[1]{Sachin Bhalekar}
		\ead{sachin.math@yahoo.co.in, sachinbhalekar@uohyd.ac.in (Corresponding Author)}

		\author[2]{Prashant M. Gade}
		\ead{prashant.m.gade@gmail.com}
		
		\address [1]{School of Mathematics and Statistics, University of Hyderabad, Hyderabad, 500046 India}
		\address [2]{Department of Physics, Rashtrasant Tukadoji Maharaj Nagpur University, Nagpur}

		\begin{abstract}
			We consider the stability of periodic map with period-$2$ in linear
			fractional difference equations where the 
			function is $f(x)=ax$ at even times and $f(x)=bx$ at odd
			times. The stability of such a map for an integer order map
			depends on product $ab$. The conditions are much
			complex for fractional maps and depend on $ab$ as
			well as $a+b$. There are no superstable period-2 orbits. 
			These conditions are useful in obtaining
			stability conditions of asymptotically periodic orbits
			with period-$2$ in the nonlinear case.
			The stability conditions are demonstrated numerically.
			The formalism can be generalized to higher periods.
		\end{abstract}
	\end{frontmatter}
	
	\section{Introduction} 
	
	The stability analysis of the fixed point is extremely helpful
	in the study of nonlinear maps of integer order. 
	The stability of periodic orbits is no less important.
	In fact, one of the definitions of
	chaos is as follows:  if $V$ is a set
	and $F: V \rightarrow  V$ is chaotic if it
	has a) sensitive dependence on initial conditions, b)
	is topologically transitive and c) the periodic points
	of $F$ are dense on $V$\cite{devaney}. 
	It has been shown that conditions b) and c) imply a).
	Thus periodic points are crucial in the 
	theory of chaos. Several 
	invariant properties of 
	the chaotic attractor can be computed using 
	unstable periodic orbits\cite{cvitan}.
	The multifractal spectrum of the attractor
	can be computed using unstable periodic
	orbits\cite{grebogi}. 
	In a very striking result,  it has been shown that
	the statistical properties of turbulence can be computed 
	using only one unstable periodic orbit\cite{kawahara,kawasaki}.
	We can compute Lyapunov exponent of the system
	using unstable periodic orbit\cite{lyapunov}.
	Several dynamic quantities can be computed using eigenvalues of a few short fundamental cycles because they are structured
	hierarchically. Longer 
	cycles only offer higher order corrections
	to these quantities\cite{cvitan}. Apart from invariant density, fractal
	dimension and Lyapunov exponents, we can
	compute the topological and metric entropy of the attractor\cite{fujisaka1983}.
	Periodic orbits form a skeleton of chaotic attractors
	and control schemes such as the Ott-Grebogi-Yorke scheme 
	have been used in controlling chaos and stabilizing
	a particular periodic orbit\cite{ogy}. Several
	bifurcations in the system, such as crisis, can be explained by 
	understanding the periodic orbits and
	their stable and unstable manifold\cite{crisis}. 
	(We note that 
	unstable periodic orbits are an
	important theoretical tool in studying quantum chaos 
	as well\cite{qc}.)   
	They have important applications, including
	control of cardiac chaos\cite{cardiac}. Methods based on 
	the detection of unstable periodic orbits have been used
	to establish low dimensional chaos in crayfish caudal
	photoreceptor \cite{crayfish}. They are used in 
	the characterization, control and prediction
	of
	experimental systems
	\cite{laser,nmr}. In short, the importance of periodic orbit in
	the theory of nonlinear dynamics and chaos cannot be overemphasized.
	
	In fractional order systems, the chaos theory is not as well developed as in integer order
	systems. However, specific results about the stability of the fixed point
	are obtained.
	We linearize around the fixed point 
	and ensure stability if the eigenvalues are inside the
	unit circle. In integer order maps,  the chain rule is applicable.
	Thus, we can 
	study the stability of fixed points of the function
	$f^n(x)$ for n-period orbits. We again linearize and the 
	stability of the n-period orbit $(x_1,x_2,\ldots, x_n)$ of
	1-d map is dictated by the condition that $\vert f'(x_n)f'(x_{n-1})
	\ldots f'(x_1)\vert <1$.
	
	Unfortunately, these conditions do not work for fractional order maps
	even for fixed points and the stability 
	is ensured if the eigenvalues are inside
	the cardioid-shaped stability region in
	the complex plane\cite{fixed}.
	This work shows that the conditions are even more
	complicated for periodic points. (In fractional 
	maps, we have only asymptotically periodic points, not strictly periodic ones).
	While it is indeed true that the fixed points are 
	given by $f(x)=x$ even for fractional order map, the
	period-n orbit is not given by roots of equation $f^n(x)=x$.
	The fixed points of the twice iterated nonlinear map $f^2(x)$ 
	do not give the 2-period orbit for the fractional order map, which is reached asymptotically.
	We will explicitly solve the system for 
	an asymptotically 2-period orbit.
	Besides, the stability  of orbit $(x_1,x_2)$
	is not given by product $f'(x_1) f'(x_2)$. It is also dependent 
	in $f'(x_1)+f'(x_2)$. Thus the stability analysis of 
	higher period orbits for fractional maps is much more
	complicated than integer order maps.
	The bifurcation diagrams can be complicated as well. 
	In integer order maps, we observe a period doubling
	cascade. When certain period becomes unstable, we obtain
	the stable solution with twice the period. In fractional
	order maps, the fixed point and period two orbit can
	both be stable at same parameter value. 
	
	In this work, we first derive the analytic conditions for the stability 
	of the periodic map. The map is linear. However, it is different for
	odd and even times. We find that the same conditions work for 
	linearized asymptotically period-two orbits of
	fractional nonlinear maps.
	
	\section{Preliminaries} \label{prel}
	In this section, we present some basic definitions and results.
	Let $ h > 0 ,\; a \in \mathbb{R}$, $ (h\mathbb{N})_a = \{ a, a+h, a+2h, \ldots\} $ and $\mathbb{N}_a=\{a, a+1, a+2, \ldots \}$.
	\begin{Def}(see \cite{ferreira2011fractional, bastos2011discrete, mozyrska2015transform}).
		For a function $x : (h\mathbb{N})_a \rightarrow  \mathbb{R}$, the forward h-difference operator if defined as 
		$$
		(\Delta_h x)(t)=\frac{x(t+h)- x(t)}{h},$$
		where t	$ \in (h\mathbb{N})_a $.
	\end{Def}
	Throughout this paper, we take $a = 0$ and $h = 1$.
	\begin{Def}\cite{mozyrska2015transform}
		For a function  $x : \mathbb{N_\circ} \rightarrow  \mathbb{R}$ the fractional sum of order $\alpha>0$ is given by
		\begin{equation}
			(\Delta^{-\alpha}x)(t) 
			= \frac{1}{\Gamma(\alpha)}\sum_{s=0}^{n}\frac{\Gamma(\alpha+n-s)}{\Gamma(n-s+1)} x(s),
		\end{equation}
		where,	$t=\alpha+n, \; n \in \mathbb{N_\circ}$.
	\end{Def}
	
	\begin{Def}\cite{mozyrska2015transform,fulai2011existence}
		Let $\mu>0$ and $m-1<\mu<m$, where $m\in \mathbb{N}$, $m=\lceil \mu \rceil$. The $\mu$th fractional Caputo-like difference is defined as
		\begin{equation}
			\Delta^\mu x(t)= \Delta^{-(m-\mu)}\left(\Delta^m x(t)\right),
		\end{equation}
		where $t\in \mathbb{N}_{m-\mu}$ and 
		\begin{equation}
			\Delta^m x(t)=\sum_{k=0}^{m}\binom{m}{k}(-1)^{m-k}x(t+k).
		\end{equation}
	\end{Def}
	\textbf{Note:} The discrete dynamical system $x(t+1)=f(x(t))$ can be written equivalently as a difference equation $\Delta x(t)=f(x(t))-x(t)$ by subtracting the term  $x(t)$ from both sides. Further, we can generalize this difference equation by replacing the operator $\Delta$ with the operator $\Delta^\alpha$, where $0<\alpha<1$. If we consider the difference equations in the higher dimensions and if the function $f$ is linear, then we get the systems of the form Equation (\ref{linCMLSys}). This is the motivation behind the term $A-I$. Of course, we can write the matrix $A'=A-I$ and use the results discussed in the literature (e.g., \cite{abu2013asymptotic}) to analyze this system.
	\begin{Def}	 \cite{mozyrska2015transform} 
		The Z-transform of a sequence $ \{y(n)\}_{n=0}^\infty $ is a complex function given by
		$Y(z)=Z[y](z)=\sum_{k=0}^{\infty} y(k) z^{-k}$
		where $z \in \mathbb{C}$ is a complex number for which the series converges absolutely.
	\end{Def}
	
	\begin{Def}\cite{mozyrska2015transform} 
		Let $\tilde{\phi}_\alpha(n)$ be a family of binomial functions defined on $\mathbb{Z}$, parametrized by $\alpha$ defined by
		\begin{eqnarray}
			\tilde{\phi}_\alpha(n) &=& \frac{\Gamma(n+\alpha-1)}{\Gamma(\alpha) \Gamma(n)}\nonumber\\
			&=& \left(
			\begin{array}{c}
				n+\alpha-1\\
				n\\
			\end{array}
			\right)
			=(-1)^n
			\left(
			\begin{array}{c}
				-\alpha\\
				n
			\end{array}
			\right).
		\end{eqnarray}
		Then
		\begin{equation*}
			Z(\tilde{\phi}_{\alpha}(t))=\frac{1}{(1-z^{-1})^{\alpha}}, \quad |z|>1.
		\end{equation*}
	\end{Def}
	
	\begin{Def}  \cite{mozyrska2015transform} 
		The convolution $\phi*x$ of the functions $\phi$ and $x$ defined on $\mathbb{N}$ is defined as
		\begin{equation*}
			\left(\phi*x\right)(n)=\sum_{s=0}^{n}\phi(n-s)x(s)=\sum_{s=0}^{n}\phi(s)x(n-s).
		\end{equation*}
		Then the Z-transform of this convolution is 
		\begin{equation}
			Z\left(\phi*x\right)(n)=\left(Z\left(\phi\right)(n)\right) \left(Z\left(x\right)(n)\right).
		\end{equation}
	\end{Def}	
	
	
	\section{The Model and Characteristic Equation}
	Let $x : \mathbb{N_\circ} \rightarrow  \mathbb{R}$ and $f$ be a map defined by 
	\begin{equation}
		f\left(x(t)\right)=
		\begin{cases}
			a x(t), & \text{if $t$ is even}\\
			b x(t), & \text{if $t$ is odd},
		\end{cases}
	\end{equation}\label{map}
	where $a$ and $b$ are real numbers.\\
	We define the fractional order discrete dynamical system using this map as
	\begin{equation}
		x(t+1)=x(0)+\sum_{j=0}^{t}\frac{\Gamma(t-j+\alpha)}{\Gamma(\alpha) \Gamma(t-j+1)}\left[f\left(x(j)\right)-x(j)\right].\label{1}
	\end{equation}
	The traditional way to obtain the characteristic equation of the systems of the form (\ref{1}) is to take Z-transform and equate the coefficient of $Z[x(t)]$ to zero. 
	The map $f$ defined in (\ref{map}) can also be written as
	\begin{equation*}
		f\left(x(t)\right)= \frac{a+b+(-1)^t (a-b)}{2}x(t).
	\end{equation*}
	Applying Z-transform to (\ref{1}), we get
	\begin{equation}
		X(z)-z x(0)=\frac{x(0)}{1-z^{-1}}+\frac{1}{(1-z^{-1})^\alpha}\left(\frac{a+b}{2}-1\right)X(z)+\frac{1}{(1-z^{-1})^\alpha}\left(\frac{a-b}{2}-1\right)X(-z). \label{ch1}
	\end{equation}
	Note that the equation (\ref{ch1}) cannot be used to find the characteristic equation because of $X(-z)$ in the last term.\\
	An elegant way to get the solution to this problem is to separate the terms $x(t)$ with even $t$ from odd values of $t$. Let us define $p(t)=x(2t)$ and $q(t)=x(2t+1)$. The system (\ref{1}) can now be written in an equivalent form as
	\begin{eqnarray}
		p(t+1)&=&x(0)+\sum_{k=0}^{t}\frac{\Gamma(2t+1-2k+\alpha)}{\Gamma(\alpha)\Gamma(2t-2k+2)}\left[(a-1)p(k)\right]\nonumber\\
		&&\,\, +\sum_{k=0}^{t}\frac{\Gamma(2t-2k+\alpha)}{\Gamma(\alpha)\Gamma(2t-2k+1)}\left[(b-1)q(k)\right],\nonumber\\
		q(t+1)&=&x(0)+\sum_{k=0}^{t}\frac{\Gamma(2t+2-2k+\alpha)}{\Gamma(\alpha)\Gamma(2t-2k+3)}\left[(a-1)p(k)\right]\nonumber\\
		&&\,\, +(a-1)p(t+1)+\sum_{k=0}^{t}\frac{\Gamma(2t+1-2k+\alpha)}{\Gamma(\alpha)\Gamma(2t-2k+2)}\left[(b-1)q(k)\right].\label{2}
	\end{eqnarray}
	If we define $\phi_1(t)=\tilde{\phi}_{\alpha}(2t)$, $\phi_2(t)=\tilde{\phi}_{\alpha}(2t+1)$ and $\phi_3(t)=\tilde{\phi}_{\alpha}(2t+2)$ then the system (\ref{2}) can be written as
	\begin{eqnarray}
		p(t+1)&=&x(0)+(a-1)\left(\phi_2*p\right)(t)+(b-1)\left(\phi_1*q\right)(t),\label{3}\\
		q(t+1)&=&x(0)+(a-1)\left(\phi_3*p\right)(t)+(a-1)p(t+1)+(b-1)\left(\phi_2*q\right)(t).\nonumber
	\end{eqnarray}
	With a few computations, we get the Z-transforms as
	\begin{eqnarray}
		Z\left(\phi_1(t)\right)&=&\frac{\left(\sqrt{z}-1\right)^{-\alpha}+\left(\sqrt{z}+1\right)^{-\alpha}}{2z^{-\alpha/2}},\nonumber\\
		Z\left(\phi_2(t)\right)&=&\frac{\left(\sqrt{z}-1\right)^{-\alpha}-\left(\sqrt{z}+1\right)^{-\alpha}}{2z^{(-\alpha-1)/2}},\nonumber\\
		Z\left(\phi_3(t)\right)&=&-z+\frac{\left(\sqrt{z}-1\right)^{-\alpha}+\left(\sqrt{z}+1\right)^{-\alpha}}{2z^{-1-\alpha/2}},\nonumber\\
		Z\left(p(t+1)\right)&=& zP(z)-zp(0),\,\, Z\left(q(t+1)\right)= zQ(z)-zq(0), \label{4}
	\end{eqnarray}
	where $Z(p(t))=P(z)$,  $Z(q(t))=Q(z)$, $p(0)=x(0)$ and $q(0)=x(1)=ax(0)$.\\
	Applying Z-transform to the system (\ref{3}) and using (\ref{4}), we get
	\begin{eqnarray}
		\left[z-(a-1)\frac{\left(\sqrt{z}-1\right)^{-\alpha}-\left(\sqrt{z}+1\right)^{-\alpha}}{2z^{(-\alpha-1)/2}}\right]P(z)&&\nonumber\\
		-(b-1)\frac{\left(\sqrt{z}-1\right)^{-\alpha}+\left(\sqrt{z}+1\right)^{-\alpha}}{2z^{-\alpha/2}}Q(z)&=& \frac{-z^2}{1-z}x(0),\nonumber\\
		(a-1)\frac{\left(\sqrt{z}-1\right)^{-\alpha}+\left(\sqrt{z}+1\right)^{-\alpha}}{2z^{-1-\alpha/2}}P(z)&&\nonumber\\
		+\left[(b-1)\frac{\left(\sqrt{z}-1\right)^{-\alpha}-\left(\sqrt{z}+1\right)^{-\alpha}}{2z^{(-\alpha-1)/2}}-z\right]Q(z)&=& \frac{z^2}{1-z}x(0).\label{5}
	\end{eqnarray}
	The characteristic equation of the system (\ref{2}) (and hence of the system (\ref{1})) can now be obtained by equating the determinant of coefficients of the terms $P(z)$ and $Q(z)$ in the system (\ref{5}) to zero as below:
	\begin{equation}
		-z(z-1)^\alpha-\frac{1}{2}(a+b-2)z^{\frac{1+\alpha}{2}}\left[\left(\sqrt{z}-1\right)^{-\alpha}-\left(\sqrt{z}+1\right)^{-\alpha}\right]+(a-1)(b-1)z^\alpha=0.\label{ch2}
	\end{equation}
	
	
	\section{Stable Region} \label{sec4}
	The zero solution of system (\ref{1}) is locally asymptotically stable if and only if all the roots $z$ of the characteristic equation (\ref{ch2}) satisfy $|z|<1$.
	Therefore, the boundary of the stable region of the system  (\ref{1}) can be obtained by substituting $z=e^{\iota t}$ in the characteristic equation (\ref{ch2}).
	We have,
	\begin{eqnarray}
		z-1&=&e^{\iota t}-1=2\sin(t/2)e^{\iota(\pi+t)/2},\nonumber\\
		\sqrt{z}-1 &=& e^{\iota t/2}-1=2\sin(t/4)e^{\iota(2\pi+t)/4},\nonumber\\
		\sqrt{z}+1 &=& e^{\iota t/2}+1=2\cos(t/4)e^{\iota t/4}. \label{6}
	\end{eqnarray}
	Using (\ref{6}), we can rewrite the characteristic equation (\ref{ch2}) as
	\begin{eqnarray}
		-2^{\alpha}\left(\sin(t/2)\right)^\alpha e^{\iota \left[t+\alpha(\pi+t)/2\right]}+(a-1)(b-1)e^{\iota \alpha t}&& \label{7}\\
		-\frac{1}{2}(a+b-2)e^{\iota t (1+\alpha)/2} 2^{\alpha}\left[\left(\sin(t/4)\right)^\alpha e^{\iota \alpha(2\pi+t)/4}-\left(\cos(t/4)\right)^\alpha  e^{\iota \alpha t/4}\right]&=&0. \nonumber 
	\end{eqnarray}
	Separating real and imaginary parts in (\ref{7}), we get
	\begin{eqnarray}
		-2^{\alpha}\left(\sin(t/2)\right)^\alpha \cos\left(\frac{\alpha (\pi+t)}{2}+t\right) + (a-1)(b-1)\cos(\alpha t)&&\nonumber\\
		- 2^{\alpha-1}(a+b-2)\left(\sin(t/4)\right)^\alpha \cos\left(\frac{\alpha \pi}{2}+t\left(\frac{1}{2}+\frac{3\alpha}{4}\right) \right)&&\nonumber\\ 
		+  2^{\alpha-1}(a+b-2)\left(\cos(t/4)\right)^\alpha \cos\left(t\left(\frac{1}{2}+\frac{3\alpha}{4}\right) \right) &=& 0,\label{8}\\
		-2^{\alpha}\left(\sin(t/2)\right)^\alpha \sin\left(\frac{\alpha (\pi+t)}{2}+t\right) + (a-1)(b-1)\sin(\alpha t)&&\nonumber\\
		- 2^{\alpha-1}(a+b-2)\left(\sin(t/4)\right)^\alpha \sin\left(\frac{\alpha \pi}{2}+t\left(\frac{1}{2}+\frac{3\alpha}{4}\right) \right)&&\nonumber\\ 
		+  2^{\alpha-1}(a+b-2)\left(\cos(t/4)\right)^\alpha \sin\left(t\left(\frac{1}{2}+\frac{3\alpha}{4}\right) \right) &=& 0.\label{9}
	\end{eqnarray}
	Equation (\ref{9}) is identically satisfied for $t=0$. For this value of $t$, the equation (\ref{8}) gives
	\begin{equation}
		b=\frac{2\left(2^\alpha-1\right)+2a\left(1-2^{\alpha-1}\right) }{2\left(a-1+2^{\alpha-1}\right)}.\label{10}
	\end{equation}
	This boundary curve (\ref{10}) can also be written as
	\begin{equation}
		\left(a-\left[1-2^{\alpha-1}\right]\right)\left(b-\left[1-2^{\alpha-1}\right]\right)=4^{\alpha-1}. \label{11}
	\end{equation}
	Let us call this boundary curve as $\Gamma_1$. The lines $a=1-2^{\alpha-1}$ and $b=1-2^{\alpha-1}$ are asymptotes for  $\Gamma_1$.\\
	If we substitute $t=\pi$, then the equations (\ref{8}) and (\ref{9}) generate the following common boundary curve $\Gamma_2$
	\begin{equation}
		b=-\frac{1+2^{\alpha}-a+2^{\alpha/2}(a-2)\sin\left(\alpha \pi/4\right) }{a-1+2^{\alpha/2}\sin\left(\alpha \pi/4\right)}.\label{12}
	\end{equation}
	Equivalently,
	\begin{equation}
		\left(a-\left[1-2^{\alpha/2}\sin\left(\alpha \pi/4\right)\right]\right)	\left(b-\left[1-2^{\alpha/2}\sin\left(\alpha \pi/4\right)\right]\right)=-2^{\alpha}\left(\cos(\alpha \pi/4)\right)^2. \label{13}
	\end{equation}
	The lines $a=1-2^{\alpha/2}\sin\left(\alpha \pi/4\right)$ and $b=1-2^{\alpha/2}\sin\left(\alpha \pi/4\right)$ are asymptotes for  $\Gamma_2$.\\
	Furthermore, the system (\ref{8})-(\ref{9}) can be solved for $a$ and $b$ as parametric functions of $t\in[0,2\pi]$. We proceed as below:\\
	Let us define $s_1=\sin^\alpha\left(t/2\right), s_2=\sin\left(t+\alpha(\pi-t)/2\right)$, $s_3=\sin^\alpha\left(t/4\right)$, $s_4=\sin\left(t(\alpha-2)/4\right)$, $s_5=\cos^\alpha\left(t/4\right)$, $s_6=\sin\left(((\alpha-2)t-2\alpha\pi)/4\right)$. Then
	\begin{eqnarray}
		a(t)&=&1+\frac{-s_1 s_2 + \sqrt{s_1\left(2^\alpha (-s_3 s_4 +s_5 s_6)(s_5 s_4-s_3 s_6)+s_1 s_2^2\right)}}{s_5 s_4- s_3 s_6},\nonumber\\
		b(t)&=&1-\frac{s_1 s_2 + \sqrt{s_1\left(2^\alpha (-s_3 s_4 +s_5 s_6)(s_5 s_4-s_3 s_6)+s_1 s_2^2\right)}}{s_5 s_4- s_3 s_6}\label{14}
	\end{eqnarray}
	is required parametric representation of the boundary curve, which we call $\Gamma_3$.
	
	\begin{The}
		For $0<\alpha<1$, the region inside the boundary curves $\Gamma_1$, $\Gamma_2$ and $\Gamma_3$ is bounded in the $ab$-plane. For any pair $(a,b)$ in this bounded region, the zero solution of system (\ref{1}) is locally asymptotically stable. 
	\end{The}
	\textbf{Proof:} 
	We assume that  $0<\alpha<1$. The intersection points between the curves $\Gamma_1$ and $\Gamma_2$ are obtained by equating the right sides of the equations (\ref{10})
	and (\ref{12}). This gives the following two points:
	\begin{eqnarray}
		(a_1,b_1)&=&\left(\frac{2^{1+\alpha/2}-2\sin(\alpha\pi/4)-\sqrt{2^\alpha +4^\alpha -2^{1+3\alpha/2}\sin(\alpha\pi/4)}}{2^{\alpha/2}-2\sin(\alpha\pi/4)},\right.\label{15}\\ &&\left.\frac{2^{1+\alpha/2}-2^{1+\alpha}\sin(\alpha\pi/4)-(2-2^\alpha)\sqrt{2^\alpha +4^\alpha -2^{1+3\alpha/2}\sin(\alpha\pi/4)}}{2^{1+\alpha/2}+2^{3\alpha/2}-2^{1+\alpha}\sin(\alpha\pi/4)-2\sqrt{2^\alpha +4^\alpha -2^{1+3\alpha/2}\sin(\alpha\pi/4)}}\right),\nonumber\\
		(a_2,b_2)&=&\left(\frac{2^{1+\alpha/2}-2\sin(\alpha\pi/4)+\sqrt{2^\alpha +4^\alpha -2^{1+3\alpha/2}\sin(\alpha\pi/4)}}{2^{\alpha/2}-2\sin(\alpha\pi/4)},\right.\label{16}\\ &&\left.\frac{2^{1+\alpha/2}-2^{1+\alpha}\sin(\alpha\pi/4)+(2-2^\alpha)\sqrt{2^\alpha +4^\alpha -2^{1+3\alpha/2}\sin(\alpha\pi/4)}}{2^{1+\alpha/2}+2^{3\alpha/2}-2^{1+\alpha}\sin(\alpha\pi/4)+2\sqrt{2^\alpha +4^\alpha -2^{1+3\alpha/2}\sin(\alpha\pi/4)}}\right).\nonumber
	\end{eqnarray}
	As the curves $\Gamma_j$ are symmetric about the line $a=b$, so are these intersection points. The point $(a_1,b_1)$ is above whereas $(a_2,b_2)$ is below the line $a=b$ and both are in the first quadrant. This also shows that there is no intersection between the curves $\Gamma_1$ and $\Gamma_2$ for the negative values of $a$ or $b$.\\
	The curves  $\Gamma_2$ and  $\Gamma_3$ intersects each other at the points $(a_3,b_3)$ and $(b_3, a_3)$, where
	\begin{eqnarray}
		a_3&=&1+2^{\alpha/2}\frac{2-\alpha+\sqrt{2(2-2\alpha+\alpha^2)\cos^2(\alpha\pi/4)}}{\alpha\cos(\alpha\pi/4)+(\alpha-2)\sin(\alpha\pi/4)},\nonumber\\
		b_3&=&1+2^{\alpha/2}\frac{2-\alpha-\sqrt{2(2-2\alpha+\alpha^2)\cos^2(\alpha\pi/4)}}{\alpha\cos(\alpha\pi/4)+(\alpha-2)\sin(\alpha\pi/4)}.
	\end{eqnarray}
	Note that $a_3<0$ and $b_3>0$.\\
	Furthermore, for any negative values of $a$ or $b$, the curve $\Gamma_3$ lies between the corresponding branches of the curves $\Gamma_1$ and $\Gamma_2$. 
	This shows that the region inside the boundary curves $\Gamma_1$, $\Gamma_2$ and $\Gamma_3$ is bounded in the $ab$-plane. Since the change in stability can occur only at these boundary curves and the system (\ref{1}) is stable at the origin, the bounded region mentioned above is the stable region for the system (\ref{1}).  This proves the result.

\begin{figure}%
		\centering
		\includegraphics[scale=1]{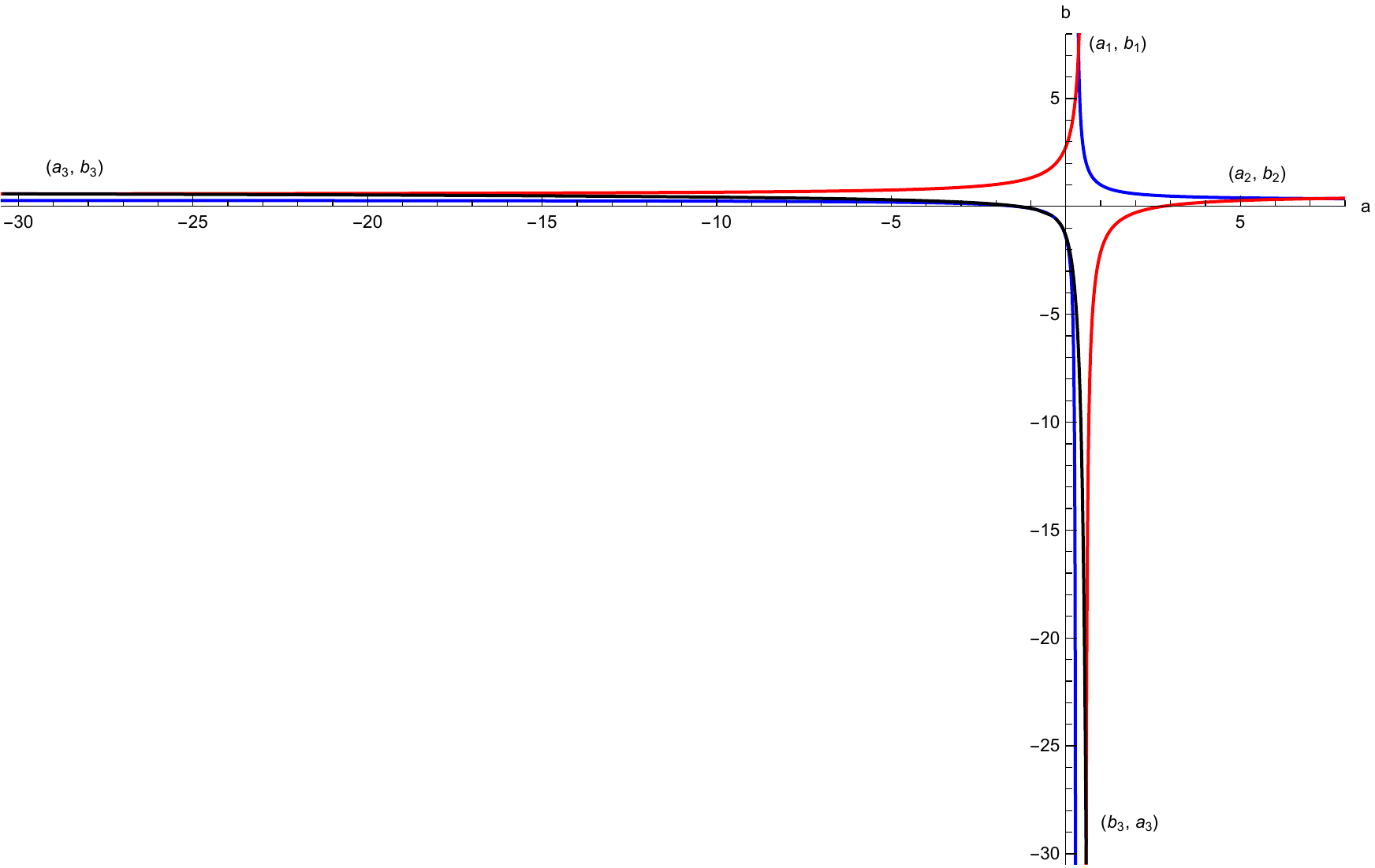} 
		\caption{The stable region of system (\ref{1}) with $\alpha=0.5$}%
		\label{f1}%
	\end{figure}

\begin{figure}%
		\centering
		\includegraphics[scale=1]{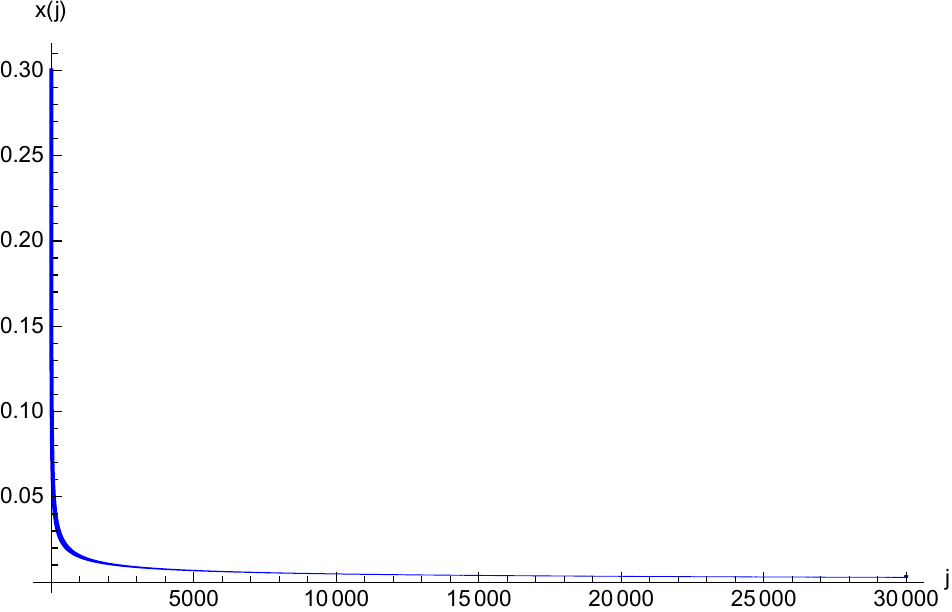} 
		\caption{The stable orbit of system (\ref{1}) with $\alpha=0.5$ and  $(a,b)=(0.6,0.7)$}%
		\label{f2}%
	\end{figure}

\begin{figure}%
		\centering
		\includegraphics[scale=1]{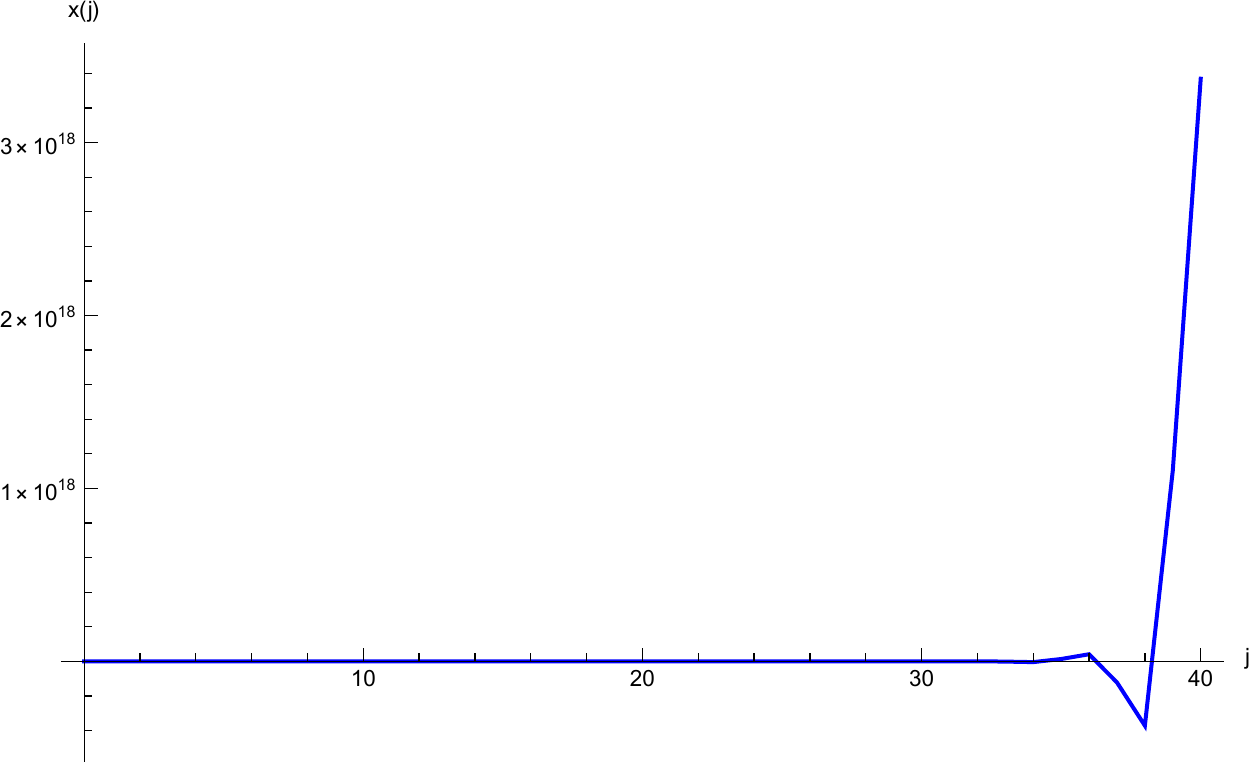} 
		\caption{The unstable orbit of system (\ref{1}) with $\alpha=0.5$ and  $(a,b)=(-2.5,3.6)$}%
		\label{f3}%
	\end{figure}

	The curves  $\Gamma_1$ (blue color), $\Gamma_2$ (red color), $\Gamma_3$ (black color) and the stable region for $\alpha=0.5$ is shown in Figure \ref{f1}.
	The stable orbit of the system  (\ref{1}) with $\alpha=0.5$ and $(a,b)=(0.6,0.7)$ inside the stable region in Figure \ref{f1} is shown in Figure \ref{f2}. On the other hand, the unstable orbit of this system with $(a,b)=(-2.5,3.6)$ outside the stable region in Figure \ref{f1} is shown in Figure \ref{f3}.
	
	Note: If $\alpha=1$, then the curve $\Gamma_3$ merges with the branch of $\Gamma_1$ in the third quadrant of $ab$-plane. Furthermore, the curves $\Gamma_1$ and $\Gamma_2$ do not intersect each other in this case and they get reduced to the curve $|ab|=1$. Note that $|ab|=1$ is the boundary of the stable region for the classical map i.e., the system $x(t+1)=f\left(x(t)\right)$, where $f$ is defined by (\ref{map}). This shows that our system and stability analysis are the continuous generalization to the classical map and the corresponding stability.

	
	\section{Application to the nonlinear systems with period-2 limit cycles}
	It is proved that \cite{kaslik2012non} the continuous-time fractional order autonomous systems of differential equations cannot have periodic solutions. However, such systems can have ``asymptotic" periodic solutions or a limit cycle \cite{kaslik2012non, bhalekar2018singular}. In this section, we show that the discrete-time fractional order systems 
	\begin{equation}
		x(t+1)=x(0)+\sum_{j=0}^{t}\frac{\Gamma(t-j+\alpha)}{\Gamma(\alpha) \Gamma(t-j+1)}\left[f\left(x(j)\right)-x(j)\right]\label{nonl}
	\end{equation}
	also have the same property. Further, we propose a necessary and sufficient condition for the existence of a period-2 limit cycle in the system  (\ref{nonl}).
	\begin{The}
		The system  (\ref{nonl}) cannot have a period-2  orbit.
	\end{The}
	Proof: If there exists the points $u$ and $v$ such that ${u,v}$ is a period-2 orbit of system (\ref{nonl}) then $x(2k)=u, x(2k+1)=v$ for $k=0,1,2,\cdots$. Therefore, for $t=0, 1$ and $2$, the system (\ref{nonl}) gives
	\begin{eqnarray}
		f(u)&=&v,\label{15a}\\
		f(v)&=&(1-\alpha)v+\alpha u,\label{16a}\\
		v&=&u+\left(\frac{\alpha(\alpha+1)}{2}+1\right)\left(f(u)-u\right)+\alpha \left(f(v)-v\right),\label{17}
	\end{eqnarray}
	respectively. Using (\ref{15a}) and (\ref{16a}) in (\ref{17}), we get
	\begin{equation}
		\left(\alpha-\frac{\alpha+1}{2}\right)\left(u-v\right)=0.
	\end{equation}
	This implies either $u=v$ or $\alpha=1$. This contradiction shows that there cannot be a period-2 orbit of the system (\ref{nonl}).
	
	\subsection{The necessary and sufficient condition for the period-2 limit cycle ${u,v}$ in the system (\ref{nonl})}
	As in (\ref{2}), we can split the system (\ref{nonl}) as
	\begin{eqnarray}
		p(t+1)&=&x(0)+\sum_{k=0}^{t}\frac{\Gamma(2t+1-2k+\alpha)}{\Gamma(\alpha)\Gamma(2t-2k+2)}\left[f(p(k))-p(k)\right]\nonumber\\
		&&\,\, +\sum_{k=0}^{t}\frac{\Gamma(2t-2k+\alpha)}{\Gamma(\alpha)\Gamma(2t-2k+1)}\left[f(q(k))-q(k)\right],\label{2a}\\
		q(t+1)&=&x(0)+\sum_{k=0}^{t}\frac{\Gamma(2t+2-2k+\alpha)}{\Gamma(\alpha)\Gamma(2t-2k+3)}\left[f(p(k))-p(k)\right]\nonumber\\
		&&\,\, +(a-1)p(t+1)+\sum_{k=0}^{t}\frac{\Gamma(2t+1-2k+\alpha)}{\Gamma(\alpha)\Gamma(2t-2k+2)}\left[f(q(k))-q(k)\right],\label{2b}
	\end{eqnarray}
	where $p(t)=x(2t)$ and $q(t)=x(2t+1)$.
	If there exists a period-2 limit cycle ${u,v}$ in the system (\ref{nonl}) then \[\lim_{t\to\infty}p(t)=u\]  and \[\lim_{t\to\infty}q(t)=v\].\\
	If $t$ is very large and $k$ is very small, then the ratios of Gamma functions in (\ref{2a}) and (\ref{2b}) become zero. On the other hand, if $k$ is very large in such cases, then $p(k)\approx u$ and $q(k)\approx v$. Therefore, taking limit as $t\to\infty$ and subtracting (\ref{2a}) from (\ref{2b}), we get
	\begin{eqnarray}
		u-v&=&(f(u)-u)\lim_{t\to\infty}\left(\sum_{k=0}^{t-1}\frac{\Gamma(2t-2k-1+\alpha)}{\Gamma(2t-2k)\Gamma(\alpha)}-\sum_{k=0}^{t}\frac{\Gamma(2t-2k+\alpha)}{\Gamma(2t-2k+1)\Gamma(\alpha)}\right)\nonumber\\
		&&+(f(v)-v) \lim_{t\to\infty}\left(\sum_{k=0}^{t-1}\frac{\Gamma(2t-2k-2+\alpha)}{\Gamma(2t-2k-1)\Gamma(\alpha)}-\sum_{k=0}^{t-1}\frac{\Gamma(2t-2k-1+\alpha)}{\Gamma(2t-2k)\Gamma(\alpha)}\right)\nonumber\\
		&=&-[(f(u)-u)-(f(v)-v)]2^{-\alpha}.\label{18}
	\end{eqnarray}
	Similarly, taking limit as $t\to\infty$ and adding (\ref{2a}) and (\ref{2b}), we get
	\begin{equation}
		u+v=2x(0)+2[(f(u)-u)+(f(v)-v)]\times[\lim_{t\to\infty}\sum_{k=0}^{2t}\frac{\Gamma(2t-k+\alpha)}{\Gamma(2t-k+1)\Gamma(\alpha)}]. \label{19}
	\end{equation}
	Since the limit in the equation (\ref{19}) tends to infinity, and all other terms are finite, we must have
	\begin{equation}
		u+v=f(u)+f(v).  \label{20}
	\end{equation}
	Solving equations (\ref{18}) and (\ref{20}), we get 
	\begin{eqnarray}
		f(u)&=&u+2^{\alpha-1}(v-u)\label{21}\\
		f(v)&=&v+2^{\alpha-1}(u-v). \label{22}
	\end{eqnarray}
	Note that for $\alpha=1$, the conditions (\ref{21})--(\ref{22}) get reduced $f(u)=v,\, f(v)=u$, the conditions for classical map $x(t+1)=f(x(t))$ to have a period-2 orbit.\\
	Recall that the period-2 orbit ${u,v}$ of the classical map is stable if $|ab|<1$, where $a=f'(u),\, b=f'(v)$. Furthermore, there are only an asymptotic period-2 orbits in the fractional order system  (\ref{nonl}) and $f(p(k))\approx f(u)+a p(k)$, $f(q(k))\approx f(v)+b q(k)$. Therefore, heuristically we can use the linearized stability analysis and expect that the point $(a,b)=(f'(u),f'(v))$ should lie inside the stable region of the system (\ref{2}).
	
	Thus the necessary and sufficient condition for the period-2 limit cycle ${u,v}$ in the system (\ref{nonl}) is \\
	(i) $f(u)=u+2^{\alpha-1}(v-u)$, $f(v)=v+2^{\alpha-1}(u-v)$, and\\
	(ii) The point $(a,b)=(f'(u),f'(v))$ lies inside the stable region bounded by the curves $\Gamma_1$, $\Gamma_2$ and $\Gamma_3$ defined in the Section \ref{sec4}.
	We verify this result with various well-known systems in the subsection below.
	
	\subsection{Examples}
	\begin{Ex}
		Consider the fractional order logistic map. In this case, we take the equation (\ref{nonl}) with $f(x)=\lambda x(1-x),$ $\lambda$ is a real parameter. 
	\end{Ex}
	The conditions (\ref{21}) and (\ref{22}) give
	\begin{eqnarray}
		(u,v)&=&\left(\frac{(2^{\alpha}-1)+\lambda + \sqrt{(\lambda-(1-2^{\alpha}))(\lambda-(1+2^{\alpha}))}}{2\lambda},\right.\nonumber \\
		&&\left. \frac{(2^{\alpha}-1)+\lambda - \sqrt{(\lambda-(1-2^{\alpha}))(\lambda-(1+2^{\alpha}))}}{2\lambda}\right). \nonumber
	\end{eqnarray}
	These points are real if $\lambda>1+2^{\alpha}$.
	
\begin{figure}%
		\centering
		\includegraphics[scale=1]{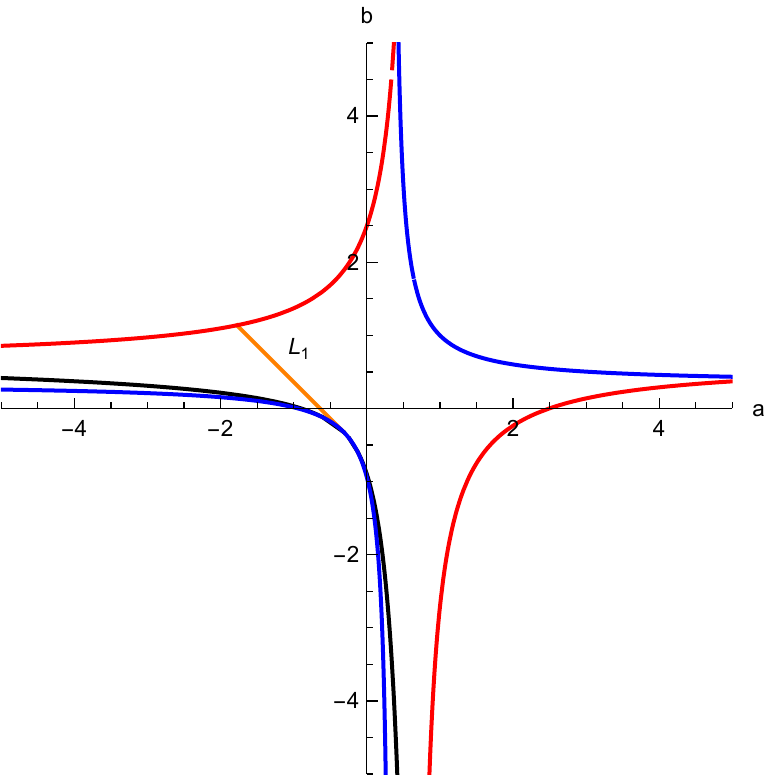} 
		\caption{The points $(a,b)=(f'(u),f'(v))$ on the line $L_1$ indicates the period-2 limit cycle in the fractional order logistic map with $\alpha=0.4$.} %
		\label{f3-1}%
	\end{figure}

\begin{figure}%
		\centering
		\includegraphics[scale=1]{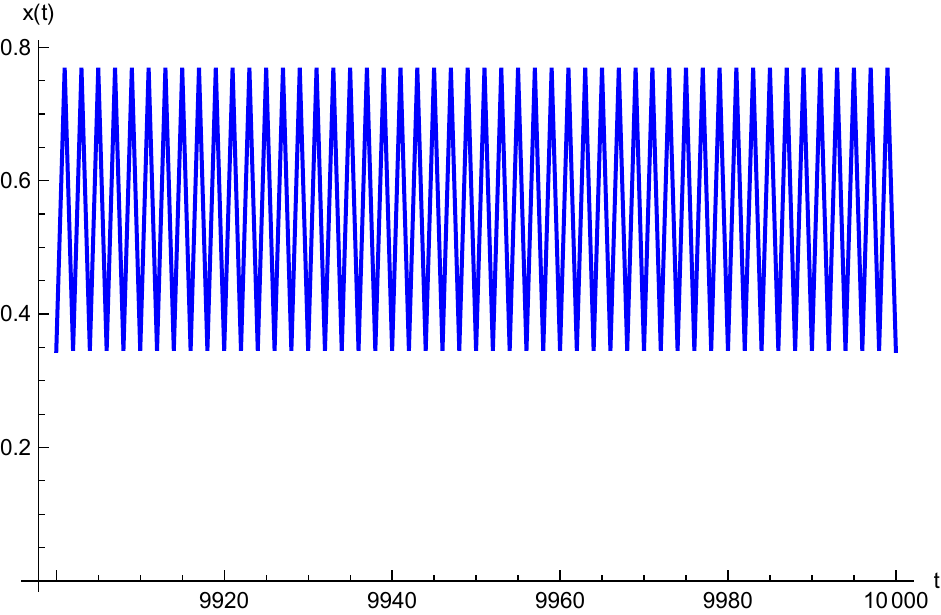} 
		\caption{Period-2 limit cycle in the fractional order  logistic map with $\alpha=0.4$ and $\lambda=2.8$.} %
		\label{f4}%
	\end{figure}

	Now, the point
	\begin{eqnarray*}
		(a,b)&=&(f'(u),f'(v))\\
		&=&\left(1-2^{\alpha}+\sqrt{(\lambda-(1-2^{\alpha}))(\lambda-(1+2^{\alpha}))},1-2^{\alpha}-\sqrt{(\lambda-(1-2^{\alpha}))(\lambda-(1+2^{\alpha}))}\right)
	\end{eqnarray*}	
	forms a straight line $L_1$, a parametric curve in $\lambda$ that intersects the curve $\Gamma_1$ at $\lambda=1+2^{\alpha}$  and $\Gamma_2$ at $\lambda=1+\sqrt{2^{\alpha}+2^{1+2\alpha}-2^{1+3\alpha/2}\sin(\alpha \pi/4)}$ (cf. Figure \ref{f3-1} for $\alpha=0.4$). 
	This shows that the fractional order logistic map has period-2 limit cycle if and only if $1+2^{\alpha}<\lambda<1+\sqrt{2^{\alpha}+2^{1+2\alpha}-2^{1+3\alpha/2}\sin(\alpha \pi/4)}$. If $\alpha=0.4$ then we need $\lambda\in(2.31951, 2.96595)$ for period-2 limit cycle. For $\lambda<2.31951$, the trajectory settles down to an equilibrium point. We can observe period-2 limit cycles when $\lambda\in(2.31951, 2.96595)$ as expected (cf. Figure \ref{f4}). Period-doubling is observed for $\lambda > 2.96595$.
	
	\begin{Ex}
		Now, we consider the fractional order cubic map. We take the equation (\ref{nonl}) with $f(x)=\beta x(6-x^2),$ where $\beta<0$ is a real parameter. 
	\end{Ex}
	We get the three expressions $(u_0, v_0)$, $(u_1, v_1)$ and $(-u_1, -v_1)$ for the points $(u,v)$ by using the conditions (\ref{21}) and (\ref{22}), where
	\begin{eqnarray*}
		u_0&=& \sqrt{(6\beta +2^{\alpha}-1)/\beta},\, v_0=-u_0,\\
		u_1 &=& \frac{1}{2}\sqrt{12+\frac{2^{\alpha}-2}{\beta}-z},\\
		v_1  &=& \frac{2^{-1-\alpha}u_1\left((2^\alpha-2)\beta+(12+z)\beta^2\right)}{\beta},\\
		z  &=& \frac{-\sqrt{4-2^{2+\alpha}-3\times 4^{\alpha}+24(2^\alpha-2)\beta+144\beta^2}}{\beta}.
	\end{eqnarray*}
	Note that, the points $(a_0, b_0)=(f'(u_0),f'(v_0))$ form a straight-line $L_2$ in the $ab-$plane, where $a_0=b_0=-12\beta-3(2^\alpha-1)$. This line $L_2$ intersects both the branches of boundary curve $\Gamma_1$ at the parameter values $\beta_0=(2-3\times 2^\alpha)/12$ and $\beta_1=(1-2^\alpha)/6$. For $\alpha=0.7$, we have  $\beta_0=-0.23946$ and $\beta_1=-0.104084$ respectively. If $\beta>\beta_1$ then the numbers $a_0$ and $b_0$ are not real; whereas if $\beta<\beta_0$ then the points $(a_0, b_0)$ are outside the stable region and hence the corresponding points $(u_0, v_0)$ don't form period-2 limit cycle for this system.\\
	
	The points $(a_1, b_1)=(f'(u_1),f'(v_1))$ form a curve $L_3$ defined by $4\left(a^2+b^2\right)+\left(9\times 2^\alpha-18\right)(a+b)+10ab+18-9\times 2^{1+\alpha}=0$ in the $ab-$plane. This curve $L_3$ intersects the boundary curves $\Gamma_1$ and $\Gamma_2$ at the parameter values $\beta_0$ and $\beta_2$, where
	\begin{eqnarray*}
		\beta_2 &=& \frac{1}{384}\left(64-33\times 2^\alpha +2^{1+\alpha/2}\sin(\alpha\pi/4)+2^{\alpha/2}\nu_1 -3\sqrt{2}\nu_2\right.\\
		&& -3\times 2^{(2+\alpha)/4}\left(179\times 2^{\alpha/2}+189\times 2^{1+3\alpha/2}-19\times 2^{\alpha/2}\cos(\alpha \pi/2)\right.\\
		&&-25\times 3^{\alpha+3}\sin(\alpha\pi/4)+2^{1+\alpha/2}\sin^2(\alpha\pi/4)+5(2^{2+\alpha/2}-2^{1+\alpha})\nu_1\\
		&&\left.\left.+3\times 2^{(1+\alpha)/2}\nu_2-3\sqrt{2}\nu_1\nu_2-6\sqrt{2}\sin(\alpha\pi/4)\nu_2\right)^{1/2}\right),\\
		\nu_1&=& \left(34+81\times 2^\alpha-2\cos(\alpha\pi/2)-9\times 2^{2+\alpha/2}\sin(\alpha\pi/4)\right)^{1/4},\\
		\nu_2 &=& 2^{\alpha/2}\left(18+9\times 2^{\alpha}-2\cos(\alpha\pi/2)-2^{\alpha/2}\nu_1+2(-5\times 2^{1+\alpha/2}+\nu_1)\sin(\alpha\pi/4)\right)^{1/2}.
	\end{eqnarray*}
	For $\alpha=0.7$, $\beta_2=-0.277584$.
	
\begin{figure}%
		\centering
		\includegraphics[scale=1]{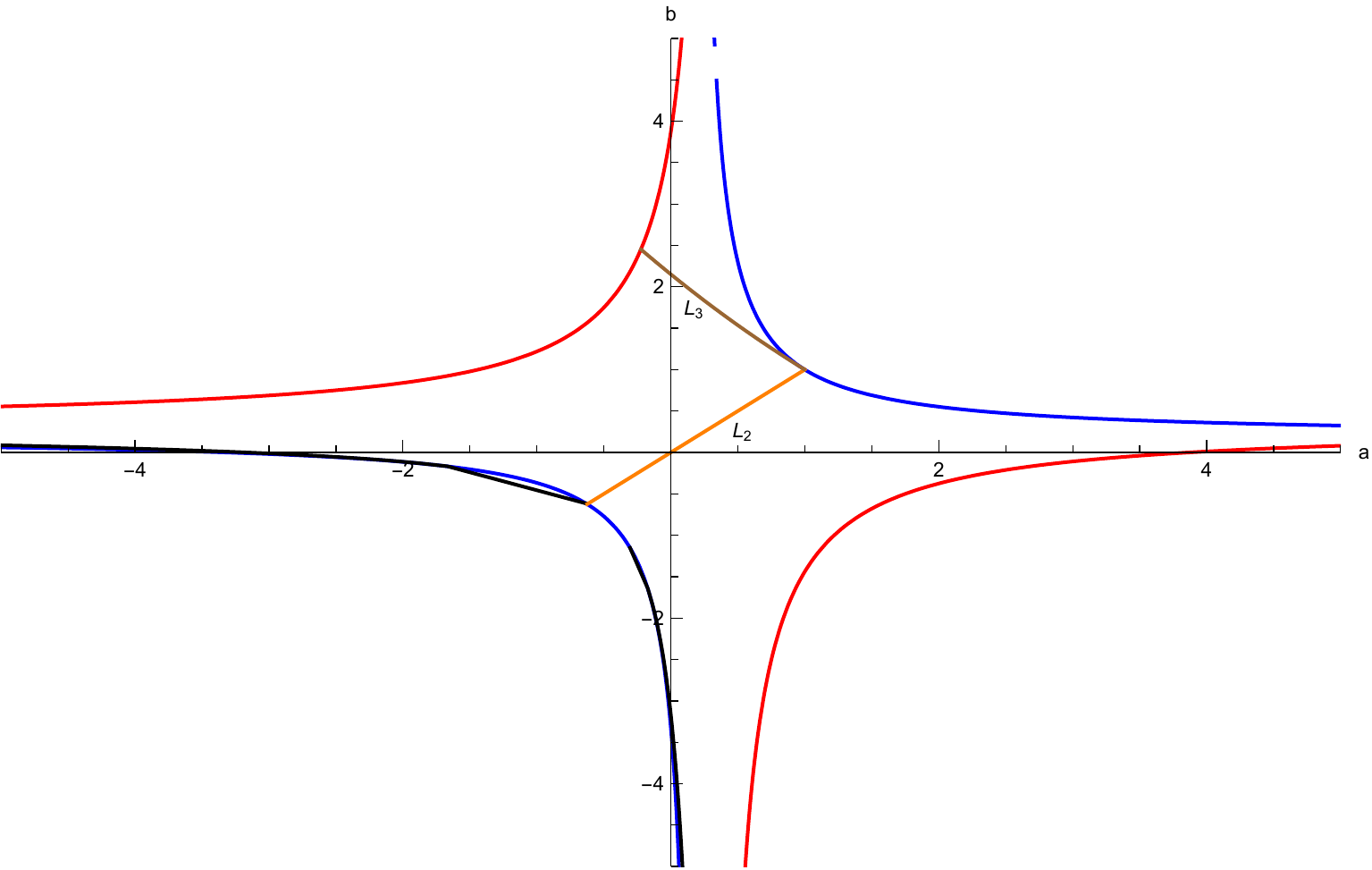} 
		\caption{The curves $L_2$ and $L_3$ in the stable region for $\alpha=0.7$ and $-0.277584<\beta<-0.104084$.} %
		\label{f5}%
	\end{figure}

	Thus, the condition for the existence of period-2 limit cycle in the fractional order cubic map is $\beta\in(\beta_2, \beta_1)$. Figure \ref{f5} shows the curves $L_2$ and $L_3$ in the stable region for $\alpha=0.7$ and $-0.277584<\beta<-0.104084$. 
	
\begin{figure}%
		\centering
		\includegraphics[scale=1]{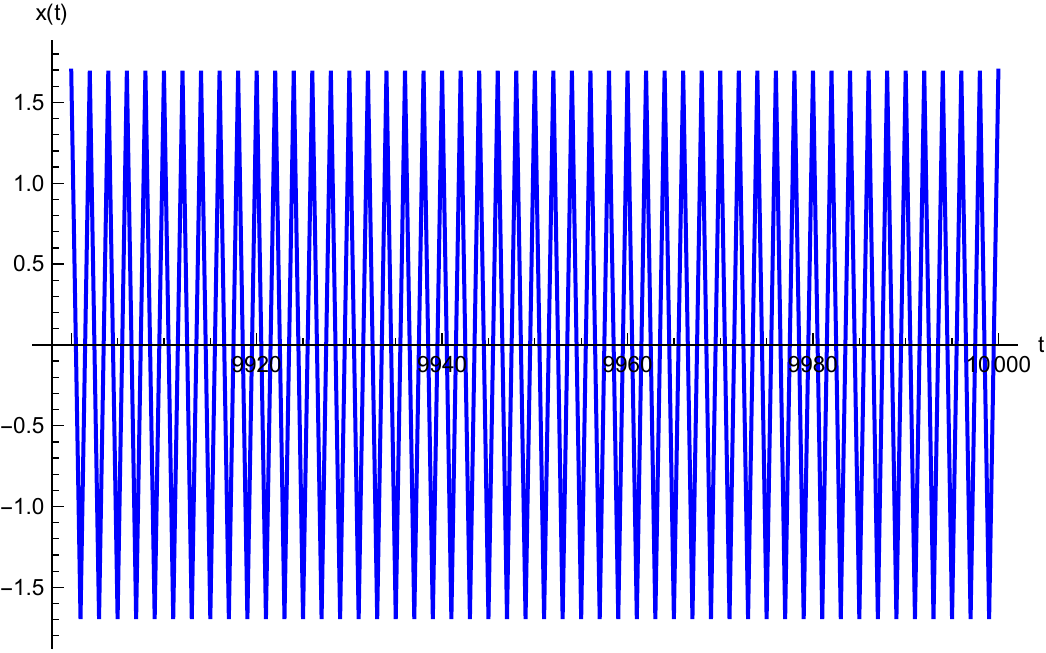} 
		\caption{Asymptotic period-2 orbit in the fractional order cubic map for $\alpha=0.7$ and $\beta=-0.20$.} %
		\label{f6}%
	\end{figure}

\begin{figure}%
		\centering
		\includegraphics[scale=1]{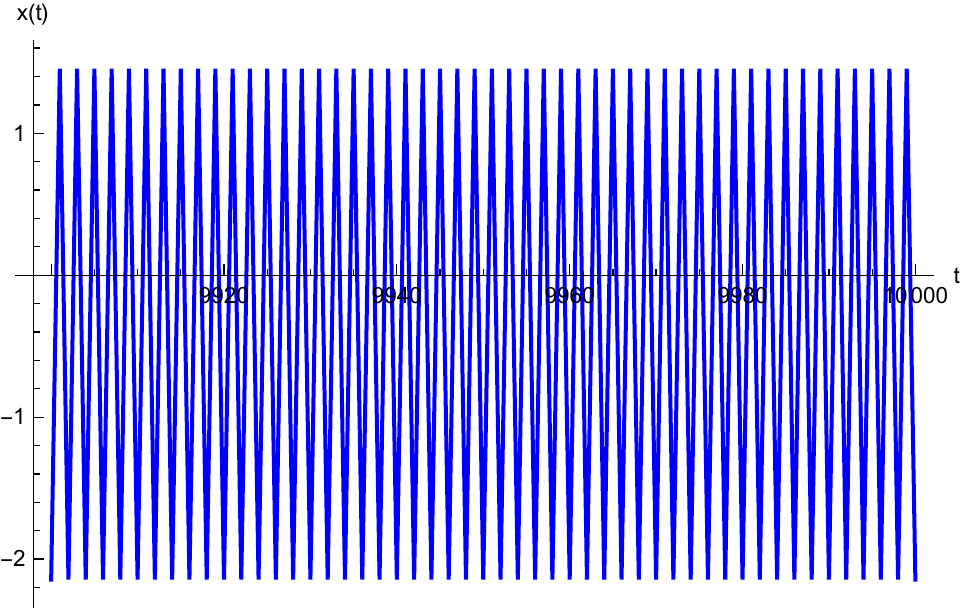} 
		\caption{Asymptotic period-2 orbit in the fractional order cubic map for $\alpha=0.7$ and $\beta=-0.26$.} %
		\label{f7}%
	\end{figure}

\begin{figure}%
		\centering
		\includegraphics[scale=1]{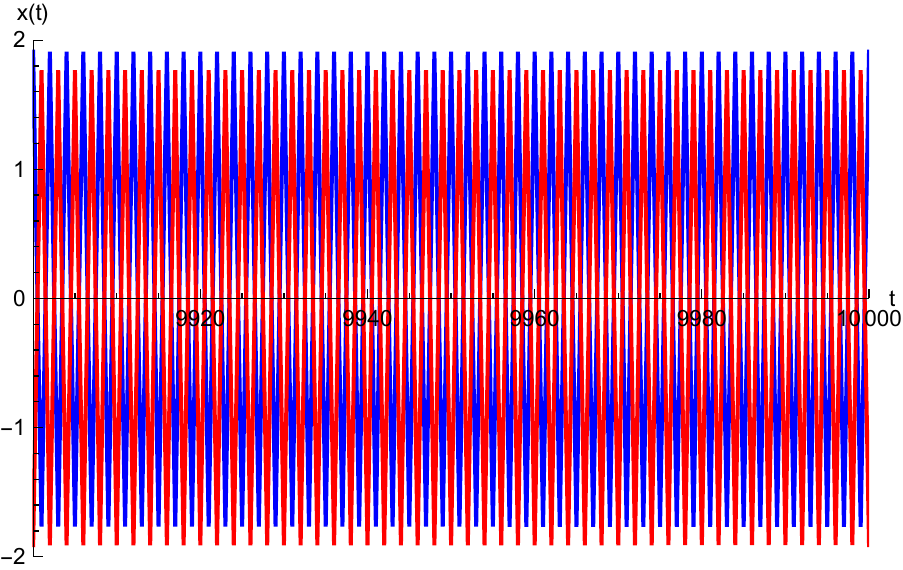} 
		\caption{Coexisting asymptotic period-2 orbits in the fractional order cubic map for $\alpha=0.7$ and $\beta=-0.24$.} %
		\label{f8}%
	\end{figure}

	For $\alpha=0.7$, the trajectory of fractional order cubic map converges to the period-2 point $(u_0, v_0)=(1.6963, -1.6963)$ when $\beta=-0.20\in(\beta_0, \beta_1)$ (cf. Figure \ref{f6}) and to the period-2 point $(u_1, v_1)=(1.45647, -2.14495)$ when $\beta=-0.26\in(\beta_2, \beta_0)$ (cf. Figure \ref{f7}). Note that the point $(-u_1, -v_1)$ indicates the existence of ``coexisting" asymptotic period-2 orbits as shown in Figure \ref{f8}.

\begin{figure}%
		\centering
		\includegraphics[scale=1]{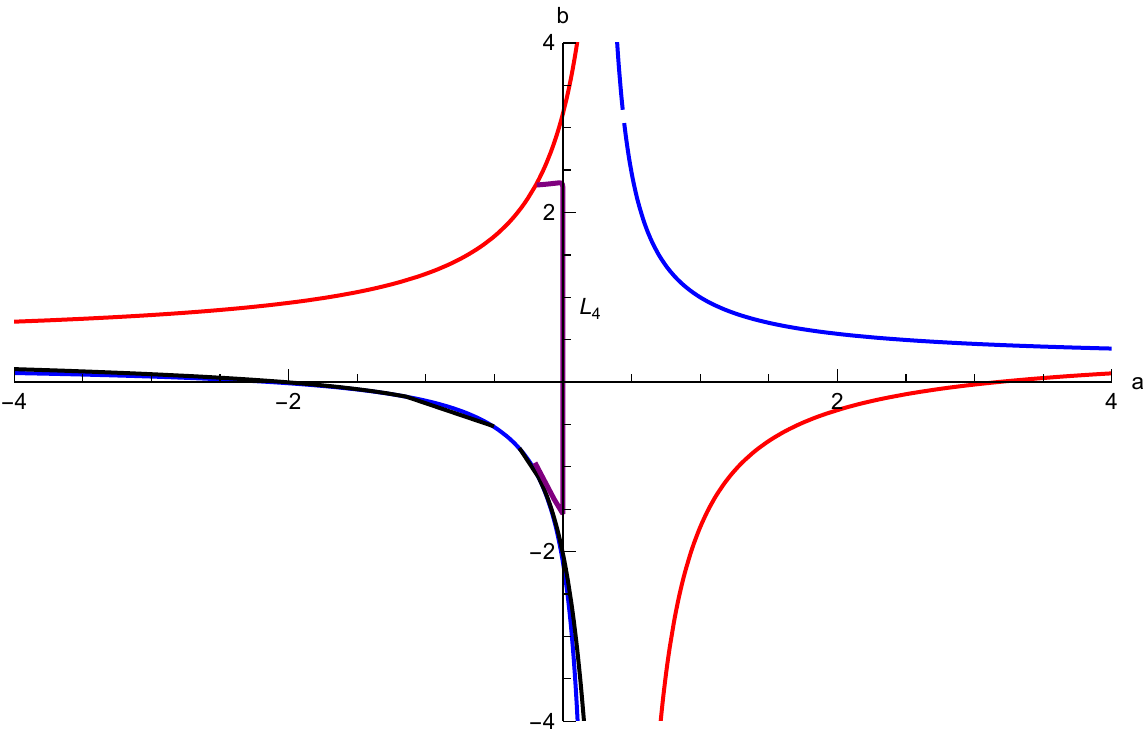} 
		\caption{The points on the curve $L_4$ for $-0.05\leq \beta \leq 0.56$ show asymptotic period-2 orbits in the fractional order Gauss map with $\alpha=0.6$.} %
		\label{f9}%
	\end{figure}

	\begin{Ex}
		In this example, we discuss the conditions for asymptotic period-2 orbits in the fractional order Gauss map (\ref{nonl}), where $f(x)=e^{-7.5x^2}+\beta$ and  $\beta$ is a real parameter.
	\end{Ex}
	Due to the transcendental nature of the function $f$, we cannot have the exact expressions for the points $(u, v)$ and $(a, b)$, unlike the previous examples. Therefore, we verify the results using numerical approximations.
	It is observed that, for $-0.05\leq \beta \leq 0.56$ and $\alpha=0.6$, the points $(a, b)=(f'(u), f'(v))$ form a curve $L_4$ that remains inside the stable region as shown in Figure \ref{f9}. The system shows asymptotic period-2 orbits for all these parameter values, as expected. We did not observe period-2 limit cycles outside this range.
	
	
	\section{Discussion}
	
	We have obtained the analytic conditions for the stability of periodic
	linear map in 
	fractional difference equations. We show that the same conditions help us
	infer the stability of asymptotically periodic orbits
	of period-$2$ in nonlinear fractional difference equations. This 
	formalism can be potentially generalized to higher periods.
	
	Unstable periodic orbits form the skeleton of chaotic attractors in 
	integer order systems. They are useful in characterization,
	prediction and control. Analysis of stable and unstable manifolds of
	periodic orbits is an indispensable tool in the theory
	of dynamical systems. The presence of chaos or the presence of stable or
	unstable manifolds of periodic orbits are open questions in
	fractional order systems. However, finding basic stability
	conditions for periodic orbit can be a useful step in
	formulating an analogous theory for fractional systems.

	\section*{Acknowledgement}
P. M. Gade thanks DST-SERB for financial assistance (Ref. CRG/2020/003993).
	
	\bibliographystyle{elsarticle-num}
	\bibliography{reference.bib}
\end{document}